\documentclass[12pt]{article}
\usepackage{amsmath,amsfonts,amsthm,anysize,mathptmx}
\marginsize{25mm}{25mm}{1cm}{1cm}
\newtheorem{theorem}{Theorem}[section]

\newtheorem{lemma}[theorem]{Lemma}
\newtheorem{proposition}[theorem]{Proposition}
\newtheorem{remark}[theorem]{Remark}
\newtheorem{corollary}[theorem]{Corollary}
\newtheorem{definition}[theorem]{Definition}

\numberwithin{equation}{section} \allowdisplaybreaks

\begin{document}

\title{\bf Existence of nontrivial solutions for  asymptotically linear periodic   Schr\"odinger
equations
 }

\author{Shaowei Chen \thanks{ E-mail address:
swchen6@163.com (Shaowei Chen)} \quad\quad Dawei Zhang
  \\ \\
\small \small\it School of Mathematical Sciences, Huaqiao
University, \\
\small Quanzhou  362021, P. R. China\\ }

\date{}
\maketitle

\begin{minipage}{13cm}
{\small {\bf Abstract:} We study the  Schr\"{o}dinger equation:
\begin{equation}
- \Delta u+V(x)u=f(x,u) ,\qquad u\in
H^{1}(\mathbb{R}^{N}),\nonumber
\end{equation} where $V$ is periodic and
$f$ is periodic in the $x$-variables, $0$ is in a gap of the
spectrum of the operator $-\Delta+V$ and $f$ is  asymptotically
linear as $|u|\rightarrow+\infty.$ We prove that under some
asymptotically linear assumptions for $f$,  this equation has a
nontrivial solution.  Our assumptions for $f$ are different from
the classical assumptions raised  by   Li and  Szulkin. \\
\medskip {\bf Key words:}  Semilinear Schr\"odinger equations;  linking; asymptotically linear.\\
\medskip 2000 Mathematics Subject Classification:  35J20, 35J60}
\end{minipage}

\section{Introduction and statement of results}\label{diyizhang}
In this paper, we consider the following Schr\"{o}dinger equation:
\begin{equation}\label{e1}
- \Delta u+V(x)u=f(x,u) ,\qquad u\in H^{1}(\mathbb{R}^{N}),
\end{equation}
where $N\geq1$, $ V (x)$ is continuous and periodic in $x_j$ for
$j = 1,\cdots,N,$ $0$ is in a gap of the spectrum of the operator
$-\Delta+V$ and $f\in C(\mathbb{R}^{N} \times\mathbb{R})$ is
 periodic in $x_j$ and   asymptotically linear  as $|u|\rightarrow\infty $.

Semilinear Schr\"odinger equations with periodic coefficients have
attracted considerable attention over the past decade. Because of
its   natural variational structure (see (\ref{vcrdfd}) in Section
\ref{nvb6ftrf} of this paper), critical point theory  is the main
method  obtaining solutions to Eq.(\ref{e1}). When $V$ is bounded
 below by a positive constant, the operator $-\Delta+V$ is positive
 definite.
In this case,   classical theorems in critical point theory, such
as the mountain pass theorem (see, for example, \cite{Willem}),
 can be used to obtain solutions to Eq.(\ref{e1}) (see the classical paper
 \cite{CR, jean} and the more recent paper  \cite{yongqing}). However,  when $0$ is in a gap of the spectrum of
the operator $-\Delta+V$, this operator has an infinitely
dimensional negative space, and the classical linking theorems
($e.g.$, \cite{Willem}) can not be applied. To overcome this
difficulty, some new infinite-dimensional linking theorems
 were developed  (see \cite{yanheng, KS, Zou, szulkin}).
      Using these
generalized linking theorems, many results on the existence and
multiplicity of  nontrivial solutions for (\ref{e1}) have been
obtained (see \cite{liu, KS1, schechter, SW, yang, YCD}). In
\cite{KS}, Kryszewski and Szulkin proved that (\ref{e1}) has
 a nontrivial solution if $f$  satisfies the Ambrosetti-Rabinowitz
 condition,    and
  has infinitely many solutions  if    the additional assumption that $f$ is odd holds. In
\cite{LiSzulkin}, Li and Szulkin obtained a
 nontrivial solution for (\ref{e1}) if $f$ satisfies some
 asymptotically linear assumptions, and in \cite{yanheng}, Ding
 proved that if $f$ is odd, then, under the same    assumptions as   in \cite{LiSzulkin},
  (\ref{e1}) has infinitely many geometrically different
  solutions. In \cite{Zou} (see also \cite{schechter}), Schechter and Zou  combined  a generalized
  linking theorem with the monotonicity methods  of Jeanjean (see
  \cite{jean}). They obtained a nontrivial solution of (\ref{e1})
  when $f$ exhibts the critical growth. A similar  approach was applied by
  Szulkin and Zou to obtain homoclinic orbits of asymptotically
  linear Hamiltonian systems (see \cite{szulkin}).
  Finally, we
should point out  that, although these generalized linking
theorems have achieved great success in strongly indefinite
problems, there are other approaches that can  be used  to deal
with (\ref{e1}) effectively. For instance,   see \cite{Ackermann,
alama, chen, HKS, pankov, PMilan, T} and the references therein.

In \cite{LiSzulkin}, Li and Szulkin studied Eq.(\ref{e1}) under
the following assumptions:
\begin{description}
\item{$(\bf{v}).$}
 $V\in C( \mathbb{R}^{N}) $ is 1-periodic in $x_j$ for $j = 1,\cdots,N$ and 0 is in a spectral gap $(\mu_{-1}, \mu_1)$ of  $-\Delta+V$.
Denote $$\mu_0 := \min\{-\mu_{-1}, \mu_1\}.$$

\item{$(\bf{f_1}).$} $f\in C( \mathbb{R}^{N}\times\mathbb{R}) $ is
1-periodic in $x_j$ for $j = 1,\cdots,N$ and $f( x,t)/t
\rightarrow 0 $ as $t\rightarrow0$ uniformly in
$x\in\mathbb{R}^{N}$.

\item{$(\bf{f_2}).$} $f( x,t) =V_\infty(x)t+f_\infty(x,t) $, where
$V_\infty$ and $f_\infty$ are 1-periodic in $x_j$ for $j =
1,\cdots,N$, $$f_\infty(x,t)/t\rightarrow 0\quad\mbox{ uniformly
in}\quad x\in\mathbb{R}^{N}\ \mbox{as} \ |t|\rightarrow\infty,$$
and $V_\infty(x) \geq\mu$ for all $x$ and some $\mu > \mu_1$.

\item{$(\bf{f_3}).$}
$\widetilde{F}(x,t):=\frac{1}{2}tf(x,t)-F(x,t)\geq 0$ for all
$(x,t)\in\mathbb{R}^N\times\mathbb{R}$, where
$F(x,t)=\int^{t}_{0}f(x,s)ds$.

\item{$(\bf{f_4}).$} There exists $\delta\in (0, \mu_0)$ such that
if $ f(x,t)/t\geq\mu_0 -\delta$, then
$\widetilde{F}(x,t)\geq\delta$.
\end{description}

Under  assumption $\bf (f_1)$,  the zero function $u=0$ is
obviously a trivial solution of (\ref{e1}). Therefore we focus on
 finding nontrivial solutions, namely solutions $u$
of (\ref{e1}) such that $u\not \equiv 0$ in $\mathbb{R}^{N}$. In
\cite{LiSzulkin}, Li and Szulkin  obtained a nontrivial solution
of equation (\ref{e1}) under the above assumptions by applying the
generalized linking theorem (see  \cite{KS} or \cite[Chapter
6]{Willem}). After \cite{LiSzulkin}, conditions similar to $\bf
(f_4)$ have become classical assumptions for strongly indefinite
problems with asymptotically linear nonlinearities (see, for
example, \cite{yanheng} and \cite{Dingbook}).

We   consider  Eq.(\ref{e1}) under  assumptions different to $\bf
(f_4)$. More precisely, we assume:

\begin{description}
\item{$(\bf{v'}).$} $0$ is not in the spectrum of the operator
\begin{eqnarray}\label{bbc66d5er}
&&T_2:L^2(\mathbb{R}^N)\rightarrow L^2(\mathbb{R}^N),\quad
u\mapsto -\Delta u+(V-V_\infty)u,\\
 &&\mbox{with domain}\quad
D(T_2):=\{u\in L^2(\mathbb{R}^N)\ |\ T_2 u\in
L^2(\mathbb{R}^N)\}.\nonumber
\end{eqnarray}

\item{$(\bf{f'_4}).$} There exist $\kappa>0$ and $\nu\in(0,\mu_0)$
such that, for every $(x,t)\in\mathbb{R}^N\times\mathbb{R}$ with
$|t|<\kappa$,
\begin{eqnarray}\label{bc77fdtdr}
|f(x,t)|\leq\nu |t| \end{eqnarray}  and for every
$(x,t)\in\mathbb{R}^N\times\mathbb{R}$ with $|t|\geq\kappa$,
\begin{eqnarray}\label{bv77fkkl} \widetilde{F}(x,t)> 0.
\end{eqnarray}

\item{$(\bf{f'_5}).$} $\widetilde{F}(x,t)> 0$ for all
$(x,t)\in\mathbb{R}^N\times(\mathbb{R}\setminus\{0\})$.
\end{description}

 Our main results are as follows:

\begin{theorem}\label{th1}
Suppose $\bf (v)$, $\bf ( v') $, $\bf ( f_{1})-\bf ( f_{3}) $, and
$\bf ( f'_{4}) $  are satisfied. Then Eq.(\ref{e1}) has a
nontrivial solution.
\end{theorem}

It is easy to verify that $\bf (f'_5)$ and the assumption that
$f( x,t)/t \rightarrow 0 $ as $t\rightarrow0$ uniformly in
$x\in\mathbb{R}^{N}$   imply $\bf (f'_4)$. Therefore, we have the
following corollary:

\begin{corollary}\label{vd5rdf}
Suppose $\bf (v)$, $\bf ( v') $, $\bf ( f_{1})$, $\bf ( f_{2}) $,
and $\bf ( f'_{5}) $  are satisfied. Then Eq.(\ref{e1}) has a
nontrivial solution.
\end{corollary}

\begin{remark}\label{bvcttdrd} There are many functions satisfying $\bf (f'_4)$  or
 $\bf (f'_5)$ that do not satisfy $\bf (f_4)$.
 An example of such a function $f$
 can be
constructed as follows: Let $b\in\mathbb{R}$ be such that
$\frac{2}{3}b\not\in\sigma(T_2)$   and $\frac{2}{3}b>\mu_1$, where
$\sigma(T_2)$ denotes the spectrum of the operator $T_2$ defined
by (\ref{bbc66d5er}). Let
$$F(x,t)=\frac{bt^2}{3}\Big(1-\frac{1}{(1+|t|)^3}\Big)\ \mbox{and}
\
f(x,t)=F'_t(x,t)=\frac{2b}{3}t\Big(1-\frac{1}{(1+|t|)^3}\Big)+\frac{bt^2\mbox{sgn}
t}{(1+|t|)^4}.$$ It is easy to verify that
$$\widetilde{F}(x,t)=\frac{b|t|^3}{2(1+|t|)^4}>0$$
for all $(x,t)\in\mathbb{R}^N\times(\mathbb{R}\setminus\{0\})$.
However, as $|t|\rightarrow+\infty,$
$f(x,t)/t\rightarrow\frac{2}{3}b>\mu_0 $ and
$\widetilde{F}(x,t)\rightarrow 0.$ Therefore, $f$ satisfies $\bf
(f'_5)$,  but  does not satisfy $\bf (f_4)$.
\end{remark}

We use  the generalized linking theorem for a class of
parameter-dependent functionals  (see \cite[Theorem 2.1]{Zou} or
Proposition \ref{b99d443} in this paper) to
 obtain a sequence of approximate solutions for (\ref{e1}). Then,
 applying
 the main theorem in \cite{HV}, we  prove that these
 approximate solutions are bounded in $L^\infty(\mathbb{R}^N)$ and $H^1(\mathbb{R}^N)$ (see Lemma \ref{bc77ftrg} and \ref{x4esdftrg}).
 These are the two most important  steps in our proof.
 Finally, using the concentration-compactness principle, we obtain
 a nontrivial solution of (\ref{e1}).

 \medskip

 \noindent{\bf Notation.} $B_r(a)$ denotes the  open ball of radius $r$ and center $a$.
For a Banach space $E,$ we denote the dual space of $E$ by $E'$,
and denote   strong and   weak convergence in $E$  by
$\rightarrow$ and $\rightharpoonup$, respectively. For $\varphi\in
C^1(E;\mathbb{R}),$ we denote the Fr\'echet derivative of
$\varphi$ at $u$ by $\varphi'(u)$.  The Gateaux derivative of
$\varphi$ is denoted by $\langle \varphi'(u), v\rangle,$ $\forall
u,v\in E.$ $L^p(\mathbb{R}^N)$ denotes the standard $L^p$ space
$(1\leq p\leq\infty)$, and $H^1(\mathbb{R}^N)$ denotes the
standard Sobolev space with norm
$||u||_{H^1}=(\int_{\mathbb{R}^N}(|\nabla u|^2+u^2)dx)^{1/2}.$  We
use $O(h)$, $o(h)$ to mean $|O(h)|\leq C|h|,$ $o(h)/|h|\rightarrow
0$ as $|h|\rightarrow 0$.

\section{Existence of approximate solutions for Eq.(\ref{e1})}\label{nvb6ftrf}
Under the  assumptions $\bf (v)$, $\bf
(f_1)$, and $\bf (f_2)$, the functional%
\begin{equation}
\Phi(u)=\frac{1}{2}\int_{\mathbb{R}^{N}}\left\vert \nabla
u\right\vert ^{2}\mathrm{d}x
+\frac{1}{2}\int_{\mathbb{R}^{N}}V(x)u^{2}\mathrm{d}x-\int
_{\mathbb{R}^{N}}F(x,u)\mathrm{d}x\label{e}%
\end{equation}
is of class $C^1$ on $X := H^1(\mathbb{R}^N )$, and the critical
points of $\Phi$ are weak solutions of (\ref{e1}).

 Assume that $\bf (v)$ holds, and let $S=-\Delta+V$ be the
self-adjoint operator acting on $L^2(\mathbb{R}^N)$ with domain
$D(S)=H^2(\mathbb{R}^N)$.  By virtue of $\bf (v)$,  we have the
orthogonal decomposition $$L^2=L^2(\mathbb{R}^N)=L^++L^-$$ such
that $S$ is negative (resp.positive) in $L^-$(resp.in $L^+$). As
in \cite[Section 2]{yanheng} (see also \cite[Chapter
6.2]{Dingbook}), let $X=D(|S|^{1/2})$ be equipped with the inner
product
$$(u,v)=(|S|^{1/2}u, |S|^{1/2}v)_{L^2}$$and norm $||u||=|||S|^{1/2}u||_{L^2}$,
 where $ (\cdot,\cdot)_{L^2}$ denotes the
inner product of $L^2$. From $\bf (v),$ $$X=H^1(\mathbb{R}^N)$$
with equivalent norms. Therefore, $X$  continuously embeds in $L^q
(\mathbb{R}^N)$ for all $2\leq q\leq2N/(N-2)$  if $N\geq 3$ and
for all $q\geq 2$  if $N=1,2$. In addition, we have the
decomposition
$$X=X^++X^-,$$ where $X^\pm=X\cap L^\pm$ is orthogonal with
respect to both $(\cdot,\cdot)_{L^2}$ and $(\cdot,\cdot)$.
Therefore, for every $u\in X$ , there is a unique decomposition
$$ u=u^++u^-,\ u^\pm\in X^\pm$$
 with $(u^+,u^-)=0$
and
$$\int_{\mathbb{R}^N}|\nabla
u|^2dx+\int_{\mathbb{R}^N}V(x)u^2dx=||u^+||^2-||u^-||^2,\ u\in
X.$$ Moreover,
\begin{eqnarray}\label{bcvttdref}
&&\mu_{-1}||u^-||^2_{L^2}\leq||u^-||^2,\quad \forall u\in X,
\end{eqnarray} and\begin{eqnarray}\label{bcvttdref2}
&&\mu_{1}||u^+||^2_{L^2}\leq||u^+||^2,\quad \forall u\in X.
\end{eqnarray}
 The functional $\Phi$ defined by (\ref{e}) can
be rewritten as
\begin{eqnarray}\label{vcrdfd}
 \Phi(u)=\frac{1}{2}(||u^+||^2-||u^-||^2)-\Psi(u),
\end{eqnarray}
where $\Psi(u)=\int_{\mathbb{R}^N}F(x,u)  \mathrm{d}x$.

Let  $\{e^\pm_k\}$ be the   total orthonormal sequence in $X^\pm.$
Let $P:X\rightarrow X^-$, $Q:X\rightarrow X^+$ be the orthogonal
projections. We define
$$|||u|||=\max\Big\{||Qu||,\sum^{\infty}_{j=1}\frac{1}{2^{k+1}}|(Pu,e^-_k)|\Big\}$$
on $X.$
 The topology
generated by $|||\cdot|||$ is denoted by $\tau$, and all
topological notation related to it will include this symbol.

\begin{definition}\label{bcv99iuytg}
Let  $\psi\in  C^1(X; \mathbb{R})$. A sequence $\{u_n\} \subset X$
is called a Cerami sequence at  level $c$ ($(C)_c$-sequence for
short) for $\psi$, if $\psi(u_n) \rightarrow c$ and $(1 +
||u_n||)||\psi'(u_n)||_{X^*}\rightarrow 0$ as
$n\rightarrow\infty.$
\end{definition}

 For $K>1$ and  $\lambda\in[1,K]$, let
\begin{eqnarray}\label{vc6drsd}
 \Phi_\lambda(u)=\frac{1}{2}\int_{\mathbb{R}^N}(|\nabla u|^2+V_+(x)u^2)dx-\lambda\Big(\frac{1}{2}\int_{\mathbb{R}^N}V_-(x)u^2dx+\Psi(u)\Big),\
 u\in X,
\end{eqnarray}
where $V_\pm(x)=\max\{\pm V(x),0\}$, $\forall x\in\mathbb{R}^N$.
It is easy to verify that a critical point $u$ of $\Phi_\lambda$
is a weak solution of
\begin{eqnarray}\label{bv66frd}
-\Delta u+V_\lambda(x)u=\lambda f(x,u),\ u\in X,
\end{eqnarray}
where $V_\lambda=V^+-\lambda V^-.$

Let $R > r > 0$ and  $u_0 \in X^+$ with $||u_0||
 = 1$. Set
$$N = \{u \in X^+ \ |\  ||u||
 = r\},\  M = \{u \in X^- \oplus \mathbb{R}^+u_0 \ |\
||u|| \leq R\}.$$ Then, $M$ is a submanifold of $X^- \oplus
\mathbb{R}^+u_0$ with boundary $\partial M.$

\begin{proposition}
[Theorem 2.1 of \cite{Zou}] \label{b99d443} Let $K>1.$ The family
of $C^1$-functional $\{H_\lambda\}$ has the form
\begin{eqnarray}\label{bcttfrd}
H_\lambda(u)=I(u)-\lambda J(u), \ u\in X,\ \lambda\in[1,K].
\end{eqnarray}
Assume
\begin{description} \item {$(a)$} $J(u)\geq 0$, $\forall u\in X$, \item
{$(b)$} $|I(u)|+J(u)\rightarrow+\infty$ as
$||u||\rightarrow+\infty$, \item {$(c)$} for all
$\lambda\in[1,K]$, $H_\lambda$ is $\tau$-sequentially upper
semi-continuous, $i.e.,$ if $|||u_n-u|||\rightarrow 0,$ then
$$\limsup_{n\rightarrow\infty}H_\lambda(u_n)\leq H_\lambda(u),$$
and $H'_\lambda$  is weakly sequentially continuous. Moreover, $
H_\lambda$ maps bounded sets to bounded sets, \item {$(d)$} there
exist $u_0 \in X^+ \setminus\{0\} $ with $||u_0||=1$, and  $R > r
> 0$ such that for all $\lambda\in[1,K]$,
$$\inf_N H_\lambda>\sup_{\partial M} H_\lambda.$$
\end{description} Then there exists $E\subset[1,K]$ such that the Lebesgue measure
of $[1,K]\setminus E$ is zero and for every  $\lambda\in E$, there
exist $c_\lambda$  and a bounded $(C)_{c_\lambda}$-sequence
 for $H_\lambda,$ where $c_\lambda$ satisfies
\begin{eqnarray}\label{bvzcx}
\sup_{M}H_\lambda\geq\sup_{\lambda\in
E}c_\lambda\geq\inf_{\lambda\in E}c_\lambda\geq\inf_N
H_\lambda.\end{eqnarray}
\end{proposition}
Using this proposition and following the same argument as the
proof of Corollary 3.4 of \cite{szulkin}, we have the following
 lemma:
\begin{lemma} \label{nxc5xrdf}
Suppose that $\bf (v)$ and $\bf (f_1)-(f_3)$  hold. Then, there
exist $K_*>1$ and  $E\subset[1,K_*]$ such that the Lebesgue
measure of $[1,K_*]\setminus E$ is zero and, for every
$\lambda\in E$, there exist $c_\lambda$ and a bounded
$(C)_{c_\lambda}$-sequence  for $\Phi_\lambda,$ where $c_\lambda$
satisfies
$$+\infty>\sup_{\lambda\in E}c_\lambda\geq\inf_{\lambda\in
E}c_\lambda>0.$$
\end{lemma}
\noindent{\bf Proof.} For $u\in X,$ let
$$I(u)=\frac{1}{2}\int_{\mathbb{R}^N}(|\nabla u|^2+V_+(x)u^2)dx$$
and $$J(u)=\frac{1}{2}\int_{\mathbb{R}^N}V_-(x)u^2dx+\Psi(u).$$
Then,  $I$ and $J$ satisfy  assumptions $(a)$ and $(b)$ in
Proposition \ref{b99d443}, and, by (\ref{vc6drsd}),
$\Phi_\lambda(u)=I(u)-\lambda J(u)$.

From (\ref{vc6drsd}),   for $u\in X,$
\begin{eqnarray}\label{bc77dtdf}
\Phi_{\lambda}(u)&=&\frac{1}{2}\int_{\mathbb{R}^{N}}(\left\vert
\nabla u\right\vert
^{2}+V(x)u^{2})dx-\frac{\lambda-1}{2}\int_{\mathbb{R}^{N}}V_-(x)u^2dx-\lambda\int
_{\mathbb{R}^{N}}F(x,u)dx\nonumber\\
&=&\frac{1}{2}||u^+||^2-\frac{1}{2}||u^-||^2-\frac{\lambda-1}{2}\int_{\mathbb{R}^{N}}V_-(x)u^2dx-\lambda\int
_{\mathbb{R}^{N}}F(x,u)dx.
\end{eqnarray}
Let $u_*\in X$ and $\{u_n\}\subset X$ be such that
$|||u_n-u_*|||\rightarrow 0$. It follows that $u_n^+\rightarrow
u^+_*$, $u^-_n\rightharpoonup u^-_*$, and $u_n\rightharpoonup
u_*$. In addition, up to a subsequence, we can assume that
$u_n\rightarrow u_*$ $a.e.$ in $\mathbb{R}^N$. Then, we have
\begin{eqnarray}
&&||u^+_n||^2\rightarrow ||u^+_*||^2,\nonumber\\
&&\liminf_{n\rightarrow\infty}||u^-_n||^2\geq
||u^-_*||^2,\nonumber\\
&&\liminf_{n\rightarrow\infty}\int_{\mathbb{R}^{N}}V_-(x)u^2_ndx\geq
\int_{\mathbb{R}^{N}}V_-(x)u^2_*dx\quad (\mbox{by the Fatou lemma}
).\nonumber
\end{eqnarray}
By the definitions of $F$ and $\widetilde{F}$, it is easy to
verify that, for all
$(x,t)\in(\mathbb{R}^N\times(\mathbb{R}\setminus\{0\})$,
$$\frac{\partial}{\partial t}\Big(\frac{F(x,t)}{t^2}\Big)=\frac{2\widetilde{F}(x,t)}{t^3}.$$
Together with $f(x,t)=o(t)$ as $|t|\rightarrow 0$ and $\bf (f_3)$,
this implies that $F(x,t)\geq 0$ for all $x$ and $t.$ By  the
Fatou lemma,
\begin{eqnarray}
\liminf_{n\rightarrow\infty}\int _{\mathbb{R}^{N}}F(x,u_n)dx\geq
\int _{\mathbb{R}^{N}}F(x,u_*)dx.\nonumber
\end{eqnarray}
Then, by  (\ref{bc77dtdf}), we obtain
\begin{eqnarray}\label{nnvb77dtdf}
\limsup_{n\rightarrow\infty}\Phi_\lambda(u_n)\leq\Phi_\lambda(u_*).\nonumber
\end{eqnarray}
This implies that $\Phi_\lambda$ is $\tau$-sequentially upper
semi-continuous.

If $u_n\rightharpoonup u_*$ in $X,$ then, for any fixed
$\varphi\in X$, as $n\rightarrow\infty,$
\begin{eqnarray}\label{mmbn88dfd}
\langle\Phi'_\lambda(u_n),\varphi\rangle
&=&\int_{\mathbb{R}^N}(\nabla u_n\nabla \varphi+V_\lambda
u_n\varphi)dx-\lambda\int_{\mathbb{R}^N}f(x,u_n)\varphi
dx\nonumber\\
&\rightarrow&\int_{\mathbb{R}^N}(\nabla u_*\nabla
\varphi+V_\lambda
u_*\varphi)dx-\lambda\int_{\mathbb{R}^N}f(x,u_*)\varphi
dx\nonumber\\
&=&\langle\Phi'_\lambda(u_*),\varphi\rangle.\nonumber
\end{eqnarray}
This implies that $\Phi'_\lambda$ is weakly sequentially
continuous. Moreover, it is easy to see that $ \Phi_\lambda$ maps
bounded sets to bounded sets. Therefore, $\Phi_\lambda$ satisfies
 assumption $(c)$  in Proposition \ref{b99d443}.

Finally, we shall verify    assumption $(d)$  in Proposition
\ref{b99d443} for $\Phi_\lambda$.

From (\ref{bc77dtdf}),  we have
\begin{eqnarray}\label{mmbn8ftfg}
\Phi_\lambda(u)=\Phi(u)-\frac{\lambda-1}{2}\int_{\mathbb{R}^{N}}V_-(x)u^2dx-(\lambda-1)\int
_{\mathbb{R}^{N}}F(x,u)dx,\quad \forall u\in X.
\end{eqnarray}
From \cite{LiSzulkin} (see also \cite[Lemma 3.1 and Lemma
3.2]{yanheng}), we know that, under assumptions $\bf (v)$ and $\bf
(f_1)-(f_3)$, there exist $u_0 \in X^+ \setminus\{0\} $ with
$||u_0||=1$, $\beta>0$, and  $R
> r
> 0$ such that
\begin{eqnarray}\label{mv88f7ft}
\inf_{N}\Phi\geq\beta\quad \mbox{and}\quad \sup_{\partial
M}\Phi\leq 0.
\end{eqnarray}
Let   $K_*>1$ be chosen such that
$$(K_*-1)\sup_{u\in N}\Big(\frac{1}{2}\int_{\mathbb{R}^{N}}V_-(x)u^2dx+\int
_{\mathbb{R}^{N}}F(x,u)dx\Big)<\beta/2.$$ Then, by
(\ref{mmbn8ftfg}) and $\inf_{N}\Phi\geq\beta$, we have that
\begin{eqnarray}\label{bbcv66dref}
\inf_{N}\Phi_\lambda\geq\beta/2,\quad \forall\lambda\in[1,K_*].
\end{eqnarray}
Moreover, by (\ref{mmbn8ftfg}) and $\sup_{\partial M}\Phi\leq 0$,
we have that
\begin{eqnarray}\label{bv77d5er11}
\sup_{\partial M}\Phi_\lambda\leq 0,\quad \forall\lambda\geq
1.\nonumber
\end{eqnarray}
Together with (\ref{bbcv66dref}), this implies that $\Phi_\lambda$
satisfies  assumption $(d)$ in Proposition \ref{b99d443} if
$\lambda\in[1,K_*]$. Therefore, for $\lambda\in[1,K_*]$,
$\Phi_\lambda$ satisfies  assumptions $(a)-(d)$ in Proposition
\ref{b99d443}. Then, the results of this lemma follow immediately
from Proposition \ref{b99d443}. \hfill$\Box$

\medskip

\begin{lemma}\label{bc6dr5df} Suppose that $\bf (v)$ and $\bf (f_1)-(f_3)$  are
satisfied.   Let $\lambda\in[1,K_*]$ be fixed, where $K_*$ is  the
constant  in Lemma \ref{nxc5xrdf}. If $\{v_{n}\}$ is a bounded
$(C)_c$ sequence for $\Phi_\lambda$ with $c\neq0$, then for every
$n\in\mathbb{N},$ there exists $a_n\in\mathbb{Z}^N$ such that, up
to a subsequence, $u_n:=v_n(\cdot+a_n)$ satisfies
\begin{eqnarray}
u_n\rightharpoonup u_\lambda\neq 0,\quad
 \Phi_\lambda(u_\lambda)\leq c \quad \mbox{and}\quad
\Phi'_\lambda(u_\lambda)= 0.
\end{eqnarray}
\end{lemma}
\noindent{\bf Proof.} The proof of this lemma is inspired by the
proof of Lemma 3.7 in \cite{szulkin}. Because $\{v_n\}$ is a
bounded sequence in $X,$   up to a subsequence, either
\begin{description}
\item {$(a)$}
$\lim_{n\rightarrow\infty}\sup_{y\in\mathbb{R}^N}\int_{B_1(y)}|v_n|^2dx=0$,
or \item {$(b)$} there exist $\varrho>0$ and $a_n\in\mathbb{Z}^N$
such that $\int_{B_1(a_n)}|v_n|^2dx\geq\varrho.$
\end{description}

If  $(a)$ occurs,  using the Lions lemma (see, for example,
\cite[Lemma 1.21]{Willem}), a similar argument as for the proof of
\cite[Lemma 3.6]{szulkin} shows that
\begin{eqnarray}\label{vvcgtdrdf}
\lim_{n\rightarrow\infty}\int_{\mathbb{R}^N}F(x,v_n)dx=0\quad
\mbox{and}\quad
\lim_{n\rightarrow\infty}\int_{\mathbb{R}^N}f(x,v_n)v^\pm_ndx=0.
\end{eqnarray}
It follows that \begin{eqnarray}\label{bcv66ftheshu}
\lim_{n\rightarrow\infty}\int_{\mathbb{R}^N}(2F(x,v_n)-f(x,v_n)v_n)dx=0.
\end{eqnarray} On the other hand, as $\{v_n\}$ is a $(C)_c-$sequence of
$\Phi_\lambda,$ we have $\langle\Phi'_\lambda(v_n),
v_n\rangle\rightarrow 0$ and $\Phi_\lambda(v_n)\rightarrow c\neq
0$. It follows that
\begin{eqnarray}
&&\int_{\mathbb{R}^N}(2F(x,v_n)-f(x,v_n)v_n)dx\nonumber\\
&=&2\Phi_\lambda(v_n)-\langle\Phi'_\lambda(v_n),
v_n\rangle\rightarrow 2c\neq0,\quad n\rightarrow\infty.
\end{eqnarray}
This contradicts (\ref{bcv66ftheshu}). Therefore,  case $(a)$
cannot occur.

If case  $(b)$ occurs, let $u_n=v_n(\cdot+a_n)$.   For every $n,$
\begin{eqnarray}\label{ggcfdrdf}
\int_{B_1(0)}|u_n|^2dx\geq\varrho.
\end{eqnarray}
Because $V$ and $F(x,t)$ are $1$-periodic in every $x_j$,
$\{u_n\}$ is still bounded in $X$,
\begin{eqnarray}\label{bbc66ftr} \lim_{n\rightarrow\infty}\Phi_\lambda(u_n)\leq
c\quad \mbox{and}\quad \Phi'_\lambda(u_n)\rightharpoonup0,\quad
n\rightarrow\infty.
\end{eqnarray}   Up to a subsequence, we  assume that
$u_n\rightharpoonup u_\lambda$ in $X$ as $n\rightarrow\infty$.
Since $u_n \rightarrow u_\lambda$ in $L^2_{ loc}( \mathbb{R}^N)$,
it follows from (\ref{ggcfdrdf}) that $u_\lambda \neq0$. Recall
that $\Phi'_\lambda(u_n)$ is weakly sequentially continuous.
Therefore,
$\Phi'_\lambda(u_n)\rightharpoonup\Phi'_\lambda(u_\lambda)$ and,
by (\ref{bbc66ftr}), $\Phi'_\lambda (u_\lambda)=0.$\hfill$\Box$

\medskip

\begin{lemma}\label{ncbvdrdf}
There exist $K_{**}>1$ and $\eta>0$ such that for any
$\lambda\in[1,K_{**}]$, if $u\neq 0$ satisfies
$\Phi'_\lambda(u)=0,$ then $||u||\geq\eta.$
\end{lemma}
\noindent{\bf Proof.} We adapt the arguments of Yang \cite[ p.
2626]{yang} and Liu \cite[Lemma 2.2]{liu}.
 Let $q=(2N-2)/(N-2)$ if $N\geq 3$ and $q=4$ if $N=1,2$.
  Note that by $\bf (f_1)$ and $\bf (f_2)$,
for any $\epsilon > 0$, there exists $C_\epsilon > 0$ such that
$$|f(x,t)|\leq \epsilon |t|+C_\epsilon|t|^{q-1}.$$
Let $u\neq 0$ be a critical point of $\Phi_\lambda$. Then, by
(\ref{bc77dtdf}) and $\langle\Phi'_\lambda(u), u^\pm\rangle=0$, we
have  that
\begin{eqnarray}\label{kvxcdxssv}
||u^\pm||^2&=&\pm(\lambda-1)\int_{\mathbb{R}^N}V_-(x)uu^\pm
dx\pm\lambda\int_{\mathbb{R}^N}f(x,u)u^\pm dx\\
&\leq&(\lambda-1)\sup_{\mathbb{R}^N}V_-\int_{\mathbb{R}^N}|u|\cdot|u^\pm|
dx\nonumber\\
&&+\epsilon\int_{\mathbb{R}^N}|u|\cdot|u^\pm|
dx+C_\epsilon\int_{\mathbb{R}^N}|u|^{q-1}|u^\pm|dx\nonumber\\
&\leq&C_1((\lambda-1)+\epsilon)||u||\cdot||u^\pm||+C_2||u||^{p-1}||u^\pm||,\nonumber
\end{eqnarray}
where $C_1$ and $C_2$ are positive constants related to the
Sobolev inequalities, and $\sup_{\mathbb{R}^N}V_-.$ From the above
two inequalities, we obtain
\begin{eqnarray}
||u||^2=||u^+||^2+||u^-||^2\leq
2C_1((\lambda-1)+\epsilon)||u||^2+2C_2||u||^{p}.
\end{eqnarray}
Because $p > 2$, this implies that $||u||
 \geq\eta$ for some $\eta> 0$ if  $\epsilon>0$ and $K_{**}-1>0$ are small enough
 and $\lambda\in[1,K_{**}]$. The desired result follows.
 \hfill$\Box$

\medskip

 Let $K=\min\{K_*,K_{**}\}$, where
$K_*$ and $K_{**}$ are the constants  that appeared  in Lemma
\ref{nxc5xrdf} and Lemma \ref{ncbvdrdf}, respectively. Combining
Lemmas $\ref{nxc5xrdf} - \ref{ncbvdrdf}$, we obtain the following
lemma:
\begin{lemma}\label{bv88f7f} Suppose $\bf (v)$  and $\bf ( f_{1})-\bf ( f_{3}) $  are satisfied. Then, there exist $\eta>0,$
 $\{\lambda_n\}\subset[1,K]$, and $\{u_n\}\subset X$ such
that $\lambda_n\rightarrow 1$,
$$\sup_{n}\Phi_{\lambda_n}(u_n)<+\infty, \quad ||u_n||\geq\eta, \quad
\mbox{and}\quad \Phi'_{\lambda_n}(u_n)=0.$$
\end{lemma}

\section{Boundedness of approximate solutions  and proofs of the main results}

In this section, we show that the sequence of approximate
solutions $\{u_n\}$ obtained in Lemma \ref{bv88f7f} is bounded in
$X$. We then give the proofs of Theorem \ref{th1} and Corollary
\ref{vd5rdf}.

\begin{lemma}\label{bc77ftrg} Suppose $\bf (v)$, $\bf ( v') $, and $\bf ( f_{1})-\bf ( f_{3}) $  are satisfied.
 Let $\{u_n\}$ be the sequence
obtained in Lemma \ref{bv88f7f}. Then, $\{u_n\}\subset
L^\infty(\mathbb{R}^N)$ and
\begin{eqnarray}\label{bc77cyfdt}
\sup_{n}||u_n||_{L^\infty(\mathbb{R}^N)}<+\infty.
\end{eqnarray}
\end{lemma}
\noindent{\bf Proof.} From $\Phi'_{\lambda_n}(u_n)=0$, we deduce
that $u_n$ is a weak solution of (\ref{bv66frd}) with
$\lambda=\lambda_n,$ $i.e.,$
\begin{eqnarray}\label{bv77ftfg}
-\Delta
u_n+V_{\lambda_n}(x)u_n=\lambda_nf(x,u_n)\quad\mbox{in}\quad
\mathbb{R}^N.
\end{eqnarray}
Because $f\in C(\mathbb{R}^N\times\mathbb{R})$ and it is
asymptotically linear, we can use  the bootstrap argument of
elliptic equations to deduce  that $u_n\in L^\infty(\mathbb{R}^N)$
and  is H\"older continuous. For every $a\in\mathbb{Z}^N,$
$u_n(\cdot+a)$ is
 still a solution of (\ref{bv77ftfg}), and so,  without loss of generality,
 we assume that for every $n\in\mathbb{N}$,  there exists $x_n\in\mathbb{R}^N$ with $|x_n|\leq
 1$ such that
 \begin{eqnarray}\label{bv77tff}
|
u_n(x_n)|=\max_{x\in\mathbb{R}^N}|u_n(x)|=||u_n||_{L^\infty(\mathbb{R}^N)}.
 \end{eqnarray}
 If (\ref{bc77cyfdt}) were not true, then
 $||u_n||_{L^\infty(\mathbb{R}^N)}\rightarrow+\infty.$  Denote
 $v_n=u_n/||u_n||_{L^\infty(\mathbb{R}^N)}.$ Then, for every $
 x\in\mathbb{R}^N$,  $|v_n(x)|\leq 1,$  and for every $n,$ $|v_n(x_n)|=1.$
 Moreover, $v_n$ satisfies
 \begin{eqnarray}\label{bvc6ftfg}
-\Delta
v_n+V_{\lambda_n}(x)v_n=\lambda_n\frac{f(x,u_n)}{u_n}v_n\quad\mbox{in}\quad
\mathbb{R}^N.
 \end{eqnarray}
As $\sup_{x\in\mathbb{R}^N,
n\in\mathbb{N}}|V_{\lambda_n}(x)|<+\infty,$ and
$\sup_{x\in\mathbb{R}^N,
n\in\mathbb{N}}\lambda_n|f(x,u_n(x))|/|u_n(x)|<+\infty$ (by $\bf
(f_1)$ and $\bf (f_3)$),  we use  the $L^p$-estimate of elliptic
equations (see, for example, \cite{GT}) to deduce that for any
$p>2$ and $R>0,$ there exists $C_R>0$ such that
\begin{eqnarray}\label{cc6tdrd} ||v_n||_{W^{2,p}(B_R(0))}\leq
C_R||v_n||_{L^p(B_{R+1}(0))}.
\end{eqnarray}
For any $x\in\mathbb{R}^N,$ $|v_n(x)|\leq 1$, which implies that
 $||v_n||_{L^p(B_{R+1}(0))}\leq|B_{R+1}(0)|^{1/p}$,
where $|A|$ denotes the Lebesgue measure of a set $A\subset
\mathbb{R}^N$. Therefore, for any $R>0$, there exists $D_R>0$ such
that
\begin{eqnarray}\label{bvttdrd}
||v_n||_{W^{2,p}(B_R(0))}\leq D_R.
\end{eqnarray}
Taking $p>N$ in (\ref{bvttdrd}), from the Sobolev embedding
theorem (see, for example, \cite[Chapter 4]{adams}), we deduce
that there exists $C'_R>0$ such that
\begin{eqnarray}\label{bc75544}
||v_n||_{C^{1,\alpha}(\overline{B_R(0)})}\leq
C'_R||v_n||_{W^{2,p}(B_R(0))}\leq C'_RD_R,
\end{eqnarray}
where $\alpha=1-N/p.$ For every $R>0,$ the embedding
$C^{1,\alpha}(\overline{B_R(0)})\hookrightarrow
C^{1}(\overline{B_R(0)})$ is compact, and so we can use the
diagonal process to deduce that there exist a subsequence
$\{v_{n_m}\}$ of $\{v_n\}$ and $v\in C^{1}(\mathbb{R}^N)$, such
that, for every $k\in\mathbb{N},$
\begin{eqnarray}\label{bv77dtre}
v_{n_m}\rightarrow v\quad \mbox{in}\quad
C^1(\overline{B_k(0)}),\quad \mbox{as}\quad m\rightarrow\infty.
\end{eqnarray}
It follows that \begin{eqnarray}\label{nvb77dtd}
v_{n_m}\rightarrow v,\ a.e.\quad \mbox{in}\quad \mathbb{R}^N,
\quad \mbox{as}\quad m\rightarrow\infty.
\end{eqnarray}
Because  $|v_{n_m}(x)|\leq 1$,  $\forall x\in\mathbb{R}^N,$
(\ref{nvb77dtd}) implies that $|v(x)|\leq 1$,  $\forall
x\in\mathbb{R}^N$. In addition, from $|v_{n_m}(x_{n_m})|=1$ and
$|x_{n_m}|\leq 1$, $m=1,2,\cdots,$ we deduce that   there exists
$x_0\in\mathbb{R}^N$ with $|x_0|\leq1$ such that, up to a
subsequence, $x_{n_m}\rightarrow x_0$ as $m\rightarrow\infty$ and
$|v(x_0)|=1.$

As the sequence $\{h_n\}$ defined by
\begin{eqnarray}
h_n(x)=\left\{
\begin{array}
[c]{ll}
f(x,u_n(x))/u_n(x) , & u_n(x)\neq 0,\\
0, &  u_n(x)=0
\end{array}
\right.\label{bvvfdre}
\end{eqnarray}
is bounded in $L^\infty(\mathbb{R}^N)$, and
$L^\infty(\mathbb{R}^N)$ is the dual space of $L^1(\mathbb{R}^N)$,
 the Banach-Alaoglu theorem (see, for example, \cite[Theorem
3.15]{rudin}) implies that, up to a subsequence,  $h_n$ converges
in the weak$^*$ topology to some function $h\in
L^\infty(\mathbb{R}^N)$, $i.e.$, for any $g\in L^1(\mathbb{R}^N)$,
$$\int_{\mathbb{R}^N}h_{n}(x)g(x) dx\rightarrow
\int_{\mathbb{R}^N}h(x)g(x) dx,\quad n\rightarrow\infty.$$ Then,
by $v_{n_m}\rightarrow v$ in $C^1_{loc}(\mathbb{R}^N)$ (see
(\ref{bv77dtre})), we have that, for any $\varphi\in C^\infty_0
(\mathbb{R}^N)$, as $m\rightarrow\infty,$
\begin{eqnarray}\label{bcv66drfd}
&&\Big|\int_{\mathbb{R}^N}h_{n_m}v_{n_m}\varphi
dx-\int_{\mathbb{R}^N}hv\varphi dx\Big|\nonumber\\
&\leq&\int_{supp \varphi}|h_{n_m}|\cdot|v_{n_m}-v|\cdot|\varphi|
dx+\Big|\int_{\mathbb{R}^N}h_{n_m}v\varphi
dx-\int_{\mathbb{R}^N}hv\varphi dx\Big|\rightarrow 0,\nonumber
\end{eqnarray}
where $supp \varphi$ denotes the support of $\varphi$.
 For any $\varphi\in C^\infty_0
(\mathbb{R}^N),$ we have
\begin{eqnarray}\label{bcvttdr}
\int_{\mathbb{R}^N}\nabla v_{n_m}\nabla\varphi
dx+\int_{\mathbb{R}^N}V_{\lambda_{n_m}}(x)v_{n_m}\varphi
dx=\lambda_{n_m}\int_{\mathbb{R}^N}h_{n_m}(x)v_{n_m}\varphi
dx.\nonumber
\end{eqnarray}
As $m\rightarrow\infty$, we have
\begin{eqnarray}
&&\lambda_{n_m}\rightarrow 1,\nonumber\\
 &&\int_{\mathbb{R}^N}\nabla
v_{n_m}\nabla\varphi dx\rightarrow
\int_{\mathbb{R}^N}\nabla v\nabla\varphi dx,\nonumber\\
&&\int_{\mathbb{R}^N}V_{\lambda_{n_m}}(x)v_{n_m}\varphi
dx\rightarrow \int_{\mathbb{R}^N}V(x)v\varphi dx,\nonumber
\end{eqnarray}
and therefore,
\begin{eqnarray}
\int_{\mathbb{R}^N}\nabla v\nabla\varphi
dx+\int_{\mathbb{R}^N}V(x)v\varphi
dx=\int_{\mathbb{R}^N}h(x)v\varphi dx,\quad \forall\varphi\in
C^\infty_0 (\mathbb{R}^N).
\end{eqnarray}
It follows that $v$ solves the linear problem
\begin{eqnarray}
-\Delta v+V(x)v=h(x)v\quad \mbox{in}\quad \mathbb{R}^N.
\end{eqnarray}
Because $v\in C^1(\mathbb{R}^N)$ and $|v(x_0)|=1,$ we can deduce
that $v\neq 0.$ Moreover, as $h\in L^\infty(\mathbb{R}^N),$ by the
regularity theorem of elliptic equations (see, for example,
\cite{GT}), we have that $v\in W^{2,2}_{loc}(\mathbb{R}^N)$. Then,
by the strong unique continuation property   (as in \cite[Theorem
6.3]{Jer}), $v(x) \neq 0$ $ a.e.$ in $\mathbb{R}^N$, which implies
$|u_{n_m}(x)| \rightarrow +\infty$, $a.e.$ in $\mathbb{R}^N$.
Hence, from $\bf (f_2)$, we have that $h_{n_m}(x)\rightarrow
V_\infty(x) $ $a.e.$ in $\mathbb{R}^N$.

We now prove that  $h=V_\infty$. It  suffices  to prove that, for
any $\varphi\in C^\infty_0(\mathbb{R}^N)$,
\begin{eqnarray}\label{vdf66dter}
\int_{\mathbb{R}^N}h\varphi dx=\int_{\mathbb{R}^N}V_\infty\varphi
dx.
\end{eqnarray}
By the Egoroff theorem (see, for example, \cite{rudin2}) and
$h_{n_m}(x)\rightarrow V_\infty(x) $ $a.e.$ in $\mathbb{R}^N$, we
deduce that, for any $\epsilon>0,$ there exists a measurable set
$E_\epsilon\subset supp \varphi$  such that $|supp
\varphi\setminus E_\epsilon|<\epsilon$ and $h_{n_m}$ converges
uniformly to $V_\infty$ on $E_\epsilon$. This implies that
\begin{eqnarray}\label{bcv77ftdf}
\lim_{m\rightarrow\infty}\int_{E_\epsilon}|h_{n_m}-V_\infty|\cdot|\varphi|dx=0\nonumber
\end{eqnarray}
and \begin{eqnarray}\label{nbv88fytr} && \sup_{m}\int_{supp
\varphi\setminus E_\epsilon}|h_{n_m}-V_\infty|\cdot |\varphi|
dx\nonumber\\
&\leq&
(\sup_{m}||h_{n_m}||_{L^\infty(\mathbb{R}^N)}+||V_\infty||_{L^\infty(\mathbb{R}^N)})\int_{supp
\varphi\setminus E_\epsilon}|\varphi|dx\leq C\epsilon,\nonumber
\end{eqnarray}
where $C>0$ is a constant independent of $m.$ Therefore,
\begin{eqnarray}\label{bcvgtdrdf}
&&\limsup_{n\rightarrow\infty}\int_{\mathbb{R}^N}|h_{n_m}-V_\infty|\cdot|\varphi|dx\nonumber\\
&\leq&\limsup_{n\rightarrow\infty}\int_{E_\epsilon}|h_{n_m}-V_\infty|\cdot|\varphi|dx\nonumber\\
&&+\limsup_{n\rightarrow\infty}\int_{supp \varphi\setminus
E_\epsilon}|h_{n_m}-V_\infty|\cdot|\varphi|dx\leq
C\epsilon.\nonumber
\end{eqnarray}
Letting $\epsilon\rightarrow 0,$ we get (\ref{vdf66dter}).
 Therefore,   $v\in
L^\infty(\mathbb{R}^N)\cap C^1(\mathbb{R}^N)$ is a nonzero
solution of the linear problem
\begin{eqnarray}\label{bcv66lk}
-\Delta u+(V(x)-V_\infty(x))u=0 \quad \mbox{in}\quad \mathbb{R}^N.
\end{eqnarray}

For $1\leq p\leq\infty,$  let $T_p$ be the operator
\begin{eqnarray}
&&T_p:L^p(\mathbb{R}^N)\rightarrow L^p(\mathbb{R}^N),\quad
u\mapsto -\Delta
u+(V-V_\infty)u,\nonumber\\
&&\mbox{with domain}\quad D(T_p):=\{u\in L^p(\mathbb{R}^N)\ |\ T_p
u\in L^p(\mathbb{R}^N)\}.\nonumber
\end{eqnarray}
Because $V-V_\infty\in L^\infty(\mathbb{R}^N)$, it was proved in
\cite{HV} that $\sigma(T_p)$, the spectrum of $ T_p$, is
independent of $p \in [1,+\infty]$. In particular, we have
$\sigma(T_2)=\sigma(T_\infty)$. Assumption $\bf (v)$ implies that
$0\not\in\sigma(T_2)$. Consequently, $0\not\in\sigma(T_\infty)$.
However, as  $v\in L^\infty(\mathbb{R}^N)$ is a nonzero solution
of (\ref{bcv66lk}), we deduce that $0\in\sigma(T_\infty)$. This
induces a contradiction. Therefore,
$\sup_{n}||u_n||_{L^\infty(\mathbb{R}^N)}<+\infty,$ which
completes the proof.\hfill$\Box$

\begin{remark}\label{hnbvvcfrd}
From  Theorem 1.2 in \cite{simon2} or Theorem C.4.2 in
\cite{simon}, we can also deduce that if (\ref{bcv66lk}) has a
nonzero solution $v\in L^\infty(\mathbb{R}^N)\cap
C^1(\mathbb{R}^N)$, then $0\in\sigma(T_2).$
\end{remark}

\begin{lemma}\label{x4esdftrg} Suppose that $\bf (v)$, $\bf ( v') $, $\bf ( f_{1})-\bf ( f_{3}) $, and
$\bf ( f'_{4}) $  are satisfied. Let $\{u_n\}$ be the sequence
obtained in Lemma \ref{bv88f7f}. Then
\begin{eqnarray}\label{bc777eedt}
0<\inf_{n}||u_n||\leq\sup_{n}||u_n||<+\infty.
\end{eqnarray}
\end{lemma}
\noindent{\bf Proof.} As $\Phi'_{\lambda_n}(u_n)=0$ and $u_n\neq
0,$   Lemma \ref{ncbvdrdf} implies that $\inf_{n}||u_n||>0.$

To prove $\sup_{n}||u_n||<+\infty$, we apply an indirect argument,
and assume by contradiction that $||u_n||\rightarrow+\infty.$

Since $\Phi'_{\lambda_n}(u_n)=0$,  by (\ref{kvxcdxssv}) and
$|f(x,u_n)|\leq C|u_n|$ for some constant $C>0$ (see $\bf (f_2)$),
we have
\begin{eqnarray}\label{ncb66ftf}
0&=&\pm||u^\pm_n||^2-(\lambda_n-1)\int_{\mathbb{R}^N}V_-(x)u_nu^\pm_ndx-\lambda_n\int_{\mathbb{R}^N}f(x,u_n)u^\pm_ndx\nonumber\\
&=&\pm||u^\pm_n||^2-\int_{\mathbb{R}^N}f(x,u_n)u^\pm_ndx+(\lambda_n-1)O(||u_n||^2).\nonumber
\end{eqnarray}
It follows that
\begin{eqnarray}\label{hf55drzm}
&&||u_n||^2-\int_{\mathbb{R}^N}f(x,u_n)(u^+_n-u^-_n)dx\nonumber\\
&=&||u^+_n||^2+||u^-_n||^2-\int_{\mathbb{R}^N}f(x,u_n)(u^+_n-u^-_n)dx=(\lambda_n-1)O(||u_n||^2).
\end{eqnarray}
Set $w_n=u_n/||u_n||$. Then by (\ref{hf55drzm}),
$$||u_n||^2\Big(1-\int_{\mathbb{R}^N}\frac{f(x,u_n)}{u_n}(w^+_n-w^-_n)w_ndx\Big)=(\lambda_n-1)O(||u_n||^2).$$
And by $\lambda_n\rightarrow 1$ as $n\rightarrow\infty,$ we have
that
\begin{eqnarray}\label{bv88dtdf}
\int_{\mathbb{R}^N}\frac{f(x,u_n)}{u_n}(w^+_n-w^-_n)w_ndx\rightarrow
1,\quad n\rightarrow\infty.
\end{eqnarray}

From Lemma \ref{bv88f7f},
$$C_0:=\sup_{n}\Phi_{\lambda_n}(u_{n})<+\infty.$$
Then, by $\Phi'_{\lambda_n}(u_n)=0$ and
$\int_{\mathbb{R}^N}\widetilde{F}(x,u_n)dx=O(||u_n||^2)$, we
obtain
\begin{eqnarray}\label{bc66drdf}
C_0&\geq&2\Phi_{\lambda_n}(u_{n})-\langle
\Phi'_{\lambda_n}(u_{n}),
u_n\rangle\nonumber\\
&=&2\lambda_n\int_{\mathbb{R}^N}\widetilde{F}(x,u_n)dx\nonumber\\
&=&(\lambda_n-1)O(||u_n||^2)+2\int_{\mathbb{R}^N}\widetilde{F}(x,u_n)dx\nonumber
\end{eqnarray}
Together with $\bf (f_3)$, this implies
\begin{eqnarray}\label{nb8fdr5er}
(\lambda_n-1)O(||u_n||^2)+C_0\geq
2\int_{\mathbb{R}^N}\widetilde{F}(x,u_n)dx\geq 2\int_{\{x\ |\
b\geq|u_n(x)|\geq\kappa\}}\widetilde{F}(x,u_n)dx
\end{eqnarray}
where $$b:=\sup_{n}||u_n||_{L^\infty(\mathbb{R}^N)}.$$ From Lemma
\ref{bc77ftrg}, we have $b<+\infty.$ As the continuous function
$\widetilde{F}$ is $1$-periodic in every $x_j$ variable,  we
deduce from (\ref{bv77fkkl}) that there exists a constant $C'>0$
such that
\begin{eqnarray}\label{bcv66fgff}
\widetilde{F}(x,t)\geq C't^2,\quad \mbox{for all}\quad
\kappa\leq|t|\leq b\quad \mbox{and}\quad  x\in\mathbb{R}^N.
\end{eqnarray}
Combining (\ref{nb8fdr5er}) and (\ref{bcv66fgff}) leads to
$$(\lambda_n-1)O(||u_n||^2)+C_0\geq 2C'\int_{\{x\ |\
b\geq|u_n(x)|\geq\kappa\}}u^2_ndx.$$ Dividing both sides of this
inequality by $||u_n||^2$ and sending $n\rightarrow\infty$, we
obtain
\begin{eqnarray}\label{bbvc77ftgd}
\lim_{n\rightarrow\infty}\int_{\{x\ |\
b\geq|u_n(x)|\geq\kappa\}}w^2_ndx=0.
\end{eqnarray}

From (\ref{bc77fdtdr}), (\ref{bcvttdref}), and (\ref{bcvttdref2}),
we have that
\begin{eqnarray}\label{nncb77dtdf}
&&\int_{\{x\ |\
|u_n(x)|<\kappa\}}\Big|\frac{f(x,u_n)}{u_n}(w^+_n-w^-_n)w_n\Big|dx\nonumber\\
&\leq&\nu\int_{\{x\ |\
|u_n(x)|<\kappa\}}|(w^+_n-w^-_n)w_n|dx\nonumber\\
&\leq&\nu\int_{\mathbb{R}^N}|(w^+_n-w^-_n)w_n|dx\nonumber\\
&\leq&\nu||w_n||^2_{L^2}\leq\frac{\nu}{\mu_0}||w_n||^2=\frac{\nu}{\mu_0}<1.
\end{eqnarray}
Because $|f(x,u_n)|\leq C|u_n|$ for some constant $C>0$ (see $\bf
(f_2)$),   (\ref{bbvc77ftgd}) gives
\begin{eqnarray}\label{nv88ihj}
&&\int_{\{x\ |\
b\geq|u_n(x)|\geq\kappa\}}\Big|\frac{f(x,u_n)}{u_n}(w^+_n-w^-_n)w_n\Big|dx\nonumber\\
&\leq&C\int_{\{x\ |\
b\geq|u_n(x)|\geq\kappa\}}|(w^+_n-w^-_n)w_n|dx\nonumber\\
&\leq& C||w^+_n-w^-_n||_{L^2}\Big(\int_{\{x\ |\
b\geq|u_n(x)|\geq\kappa\}}w^2_ndx\Big)^{1/2}\nonumber\\
&\leq&C||w_n||_{L^2}\Big(\int_{\{x\ |\
b\geq|u_n(x)|\geq\kappa\}}w^2_ndx\Big)^{1/2}\rightarrow 0,\quad
n\rightarrow\infty.
\end{eqnarray}
Combining (\ref{nncb77dtdf}) and (\ref{nv88ihj}) yields
\begin{eqnarray}\label{nnv8yydfd}
&&\limsup_{n\rightarrow\infty}\int_{\mathbb{R}^N}\Big|\frac{f(x,u_n)}{u_n}(w^+_n-w^-_n)w_n\Big|dx\nonumber\\
&\leq&\limsup_{n\rightarrow\infty}\int_{\{x\ |\
|u_n(x)|<\kappa\}}\Big|\frac{f(x,u_n)}{u_n}(w^+_n-w^-_n)w_n\Big|dx\nonumber\\
&&+\limsup_{n\rightarrow\infty}\int_{\{x\ |\
b\geq|u_n(x)|\geq\kappa\}}\Big|\frac{f(x,u_n)}{u_n}(w^+_n-w^-_n)w_n\Big|dx<1.\nonumber
\end{eqnarray}
This contradicts (\ref{bv88dtdf}). Therefore, $\{u_n\}$ is bounded
in $X.$\hfill$\Box$

\bigskip

\noindent{\bf Proof of Theorem \ref{th1}.} Let $\{u_n\}$ be the
sequence obtained in Lemma \ref{bv88f7f}.  From Lemma
\ref{x4esdftrg},  $\{u_n\}$ is bounded in $X$. Therefore,  up to a
subsequence, either
\begin{description}
\item {$(a)$}
$\lim_{n\rightarrow\infty}\sup_{y\in\mathbb{R}^N}\int_{B_1(y)}|u_n|^2dx=0$,
or \item {$(b)$} there exist $\varrho>0$ and $y_n\in\mathbb{Z}^N$
such that $\int_{B_1(y_n)}|u_n|^2dx\geq\varrho.$
\end{description}
According to (\ref{vvcgtdrdf}), if  case $(a)$ occurs,
\begin{eqnarray}
\lim_{n\rightarrow\infty}\int_{\mathbb{R}^N}f(x,u_n)u^\pm_ndx=0.\nonumber
\end{eqnarray}
Then, by (\ref{kvxcdxssv}) and $\lambda_n\rightarrow 1$, we have
\begin{eqnarray}\label{gdf55ere}
||u^\pm_n||^2&=&\pm(\lambda_n-1)\int_{\mathbb{R}^N}V_-(x)u_nu^\pm_n
dx\pm\lambda_n\int_{\mathbb{R}^N}f(x,u_n)u^\pm_n dx\nonumber\\
&\leq&
C(\lambda_n-1)||u_n||^2_{L^2}+K\Big|\int_{\mathbb{R}^N}f(x,u_n)u^\pm_n
dx\Big|\rightarrow 0.
\end{eqnarray}
This contradicts $\inf_{n}||u_n||>0$ (see (\ref{bc777eedt})).
Therefore,  case $(a)$ cannot occur. As case  $(b)$ therefore
occurs,   the proof of Lemma \ref{bc6dr5df} implies that there
exists $y_n\in\mathbb{Z}^N$ such that $w_n=u_n(\cdot+y_n)$
satisfies $w_n\rightharpoonup u_0\neq 0$. Because
$\Phi'_{\lambda_n}(u_n)=0$ (by Lemma \ref{bv88f7f}), we have
$\Phi'_{\lambda_n}(w_n)=0$. From (\ref{mmbn8ftfg}), we have that,
for any $\varphi\in X,$
\begin{eqnarray}\label{bv88ftdfd}
&&\langle\Phi'_{\lambda_n}(w_n),\varphi\rangle\nonumber\\
&=&\langle\Phi'(w_n),\varphi\rangle-(\lambda_n-1)\int_{\mathbb{R}^N}V_-(x)w_n\varphi
dx-(\lambda_n-1)\int_{\mathbb{R}^N}f(x,w_n)\varphi dx.\nonumber
\end{eqnarray}
Together with  $\Phi'_{\lambda_n}(w_n)=0$ and
$\lambda_n\rightarrow 1$, this yields
$$\langle\Phi'(w_n),\varphi\rangle\rightarrow 0,\quad \forall\varphi\in X.$$
Finally, by $w_n\rightharpoonup u_0\neq 0$ and the weakly
sequential continuity of $\Phi'$, we have that $\Phi'(u_0)=0.$
Therefore, $u_0$ is a nontrivial solution of Eq.(\ref{e1}). This
completes the proof.\hfill$\Box$

\bigskip

\noindent{\bf Proof of Corollary \ref{vd5rdf}.} Assumption   $\bf
(f'_5)$ and  the assumption that $f(x,t)/t\rightarrow 0$ uniformly
in $x\in\mathbb{R}^N$  as $t\rightarrow 0$
 imply $\bf (f'_4)$. Thus, this corollary follows immediately
from Theorem \ref{th1}.\hfill$\Box$

\end{document}